\input amstex
\input amsppt.sty
\magnification=\magstep1
\hsize=32truecc
\vsize=22.5truecm
\baselineskip=16truept
\NoBlackBoxes
\TagsOnRight \pageno=1 \nologo
\def\Z{\Bbb Z}

\def\l{\left}
\def\r{\right}
\def\bg{\bigg}
\def\({\bg(}
\def\[{\bg\lfloor}
\def\){\bg)}
\def\]{\bg\rfloor}
\def\t{\text}
\def\f{\frac}

\def\ls{\leqslant}
\def\gs{\geqslant}

\def\Proof{\noindent{\it Proof}}

\def\Remark{\medskip\noindent{\it  Remark}}

\hbox {Preprint, {\tt arXiv:1209.3729}}
\bigskip
\topmatter
\title Some new inequalities for primes\endtitle
\author Zhi-Wei Sun\endauthor
\leftheadtext{Zhi-Wei Sun}
 \rightheadtext{Some new inequalities for primes}
\affil Department of Mathematics, Nanjing University\\
 Nanjing 210093, People's Republic of China
  \\  zwsun\@nju.edu.cn
  \\ {\tt http://math.nju.edu.cn/$\sim$zwsun}
\endaffil
\abstract For $n=1,2,3,\ldots$ let $p_n$ be the $n$th prime. We show
that $p_n>n+\sum_{k=1}^n p_k/k$ for all $n\gs 125$, and
$\sum_{k=1}^n kp_k<n^2p_n/3$ for all $n>30$.
\endabstract
\thanks 2010 {\it Mathematics Subject Classification}.\,Primary 11A41;
Secondary 11B75.
\newline\indent {\it Keywords}. Primes, inequalities.
\newline\indent Supported by the National Natural Science
Foundation (grant 11171140) of China and the PAPD of Jiangsu Higher
Education Institutions.
\endthanks
\endtopmatter
\document

\heading{1. Introduction}\endheading

For $n\in\Z^+=\{1,2,3,\ldots\}$ let $p_n$ denote the $n$th prime.
Robert Mandl ever conjectured that
$$\sum_{k=1}^np_k<\f{np_n}2\quad\t{for all}\ n\gs9.$$
In 1975 J. B. Rosser and L. Schoenfeld [RS] claimed to have a proof of Mandl's inequality
but they have never published the details. In 1998 P. Dusart [D98] gave a detailed proof of
the Mandl inequality. In 2007, M. Hassani [H] refined Mandl's inequality as
$$\sum_{k=1}^np_k<\f{np_n}2-\f{n^2}{14}\quad\t{for}\ n\gs10.$$

In this paper we mainly establish the following new theorem.

\proclaim{Theorem 1.1} We have the inequality
$$p_n>n+\sum_{k=1}^n\f{p_k}k\quad\t{for all}\ n\gs125.\tag1.1$$
Also,
$$\sum_{k=1}^n kp_k<\f{n^2}3p_n-\f{n^3}9\quad\t{for all}\ n\gs2700\tag1.2$$
and
$$\sum_{k=1}^np_k<\f{n}2p_n-\f{n^2}4\quad\t{for all}\ n\gs417.\tag1.3$$
\endproclaim

We also have the following general result.

\proclaim{Theorem 1.2} Let $b>0$. For any $a>-b$, there is a positive integer $n_0(a,b)$ such that
$$\sum_{k=1}^nk^{a-1}p_k^b<\f{n^{a}p_n^b}{a+b}-\f{bn^{a+1}p_n^{b-1}}{(a+b)^2}$$
for all $n\gs n_0(a,b)$.
When $n$ is sufficiently large, we also have
$$\sum_{k=1}^n\f{p_k^b}{k^{b+1}}<\f{p_n^{b+1}-np_n^b}{(b+1)n^{b+1}}.$$
\endproclaim
\Remark. We can take
$$\gather n_0(3,1)=6048,\ n_0(4,1)=6077,\ n_0(5,1)=9260,\ n_0(6,1)=20477,
\\ n_0(7,1)=30398,\ n_0(8,1)=37358,\ n_0(9,1)=37374,\ n_0(10,1)=92608.
\endgather$$
Also, we may take $n_0(-2047,1)=2215$ and $n_0(-1,2)=348271$.

\medskip

We omit the proof of Theorem 1.2 since it is quite similar to the proof of (1.1).

To conclude this section, we pose a conjecture.

\proclaim{Conjecture 1.1} {\rm (i)} If $q_1,q_2,\ldots,q_n$ are the first $n$ primes $p$ with $M_p=2^p-1$ prime, then
$\sum_{k=1}^nq_k<nq_n/\log n$.

{\rm (ii)} Let $(F_n)_{n\gs0}$ be the Fibonacci sequence. If $q(1),q(2),\ldots,q(n)$ are the first $n$ primes $p$ with $F_p$ prime, then
$\sum_{k=1}^nq(k)<n\,q(n)/\log n$.
\endproclaim
\Remark. We have verified part (i) for all known Mersenne primes and part (ii) for all known Fibonacci primes.

\heading{2. Proof of Theorem 1.1}\endheading

We first give a lemma obtained by P. Dusart [D99].

\proclaim{Lemma 2.1} We have
$$p_k>k\l(\log k+\log\log k-1+\f{\log\log k-2.25}{\log k}\r)\ \ \t{for all}\ k\gs2,\tag2.1$$
and
$$p_k\ls k\l(\log k+\log\log k-1+\f{\log\log k-1.8}{\log k}\r)\ \ \t{for all}\ k\gs27076.\tag2.2$$
\endproclaim

\proclaim{Lemma 2.2} Let $m$ and $n$ be positive integers with $3\ls m<\sqrt n$. Then
$$\sum_{k=m}^n\f1{\log k}<\f n{\log n}+\f{4n}{\log^2 n}+\f{\sqrt n-2(m-1)}{\log^2 m}.\tag2.3$$
\endproclaim
\Proof. Since
$$\align&\sum_{k=m}^n\f1{\log k}=\sum_{k=m}^n\f{k-(k-1)}{\log k}
\\=&\sum_{k=m}^{n-1}k\l(\f1{\log k}-\f1{\log(k+1)}\r)+\f n{\log n}-\f{m-1}{\log m}
\\=&\sum_{k=m}^{n-1}\f{k\log(1+1/k)}{(\log k)\log(k+1)}+\f n{\log n}-\f{m-1}{\log m}
\\<&\sum_{k=m}^{n-1}\f1{\log^2k}+\f n{\log n}-\f{m-1}{\log m}
\endalign$$
and
$$\align\sum_{k=m}^{n-1}\f1{\log^2k}=&\sum_{k=m}^{\lfloor\sqrt n\rfloor}\f1{\log^2k}+\sum_{k=\lfloor\sqrt n\rfloor+1}^{n-1}\f1{\log^2k}
\\\ls&\f{\lfloor\sqrt{n}\rfloor-(m-1)}{\log^2m}+\f{n-1-\lfloor\sqrt n\rfloor}{\log^2\sqrt n}
\\<&\f{\sqrt n-m+1}{\log^2m}+\f {4(n-\sqrt n)}{\log^2 n},
\endalign$$
we have
$$\sum_{k=m}^n\f1{\log k}<\f n{\log n}-\f{m-1}{\log m}+\f{\sqrt n-m+1}{\log^2 m}+\f {4n}{\log^2 n}$$
and hence the desired (2.3) follows. \qed

\medskip
\noindent {\it Proof of Theorem 1.1}. For $125\ls n<50000$ we can verify the desired inequality directly.
Below we assume that $n\gs 50000$.

Let $m=27076$.
Via computer we find that
$$\sum_{k=1}^{m-1}\f{p_k}k<283452.35.$$
Therefore
$$\sum_{k=1}^n\f{p_k}k-283452.35<\sum_{k=m}^n\f{p_k}k\ls\sum_{k=m}^n\l(\log k+\log\log k-1+\f{\log\log k-1.8}{\log k}\r).\tag2.4$$
by Lemma 2.1. Clearly
$$\align&\sum_{k=m}^{n-1}(\log k+\log\log k)
\\\ls&\sum_{k=m}^{n-1}\int_k^{k+1}(\log x+\log\log x)dx=\int_m^n(\log x+\log\log x)dx
\\=&x(\log x+\log\log x)\big|_{x=m}^n-\int_m^nx\l(\f1x+\f1{x\log x}\r)dx
\\=&n(\log n+\log\log n-1)-m(\log m+\log\log m-1)-\int_m^n\f{dx}{\log x}
\endalign$$
and
$$\int_m^n\f{dx}{\log x}=\f x{\log x}\bigg|_{x=m}^n-\int_m^n\l(-\f x{x\log^2x}\r)dx\gs\f{n}{\log n}-\f m{\log m}+\f{n-m}{\log^2n}.$$
Therefore
$$\aligned\sum_{k=m}^{n}(\log k+\log\log k-1)
\ls&\log n+\log\log n-(n-m+1)
\\&+n(\log n+\log\log n-1)-m(\log m+\log\log m-1)
\\&-\f n{\log n}+\f m{\log m}-\f{n-m}{\log^2 n}
\endaligned\tag2.5$$
As $n\gs 27076>164^2$, by Lemma 2.2 we have
$$\sum_{k=164}^n\f1{\log k}<\f n{\log n}+\f{4n}{\log^2 n}+\f{\sqrt n-2\times 163}{\log^2 164}$$
and hence
$$\align\sum_{k=m}^n\f1{\log k}=&\sum_{k=164}^n\f1{\log k}-\sum_{k=164}^{m-1}\f1{\log k}
\\<&\f n{\log n}+\f{4n}{\log^2 n}+\f{\sqrt n-326}{\log^2 164}-2948.64
\\<&\f n{\log n}+\f{4n}{\log^2 n}+\f{\sqrt n}{\log^2 164}-2961.17.
\endalign$$
Thus
$$\sum_{k=m}^n\f{\log\log k-1.8}{\log k}\ls(\log\log n-1.8)\l(\f n{\log n}+\f{4n}{\log^2 n}+\f{\sqrt n}{\log^2 164}-2961.17\r).\tag2.6$$
Combining (2.4)-(2.6) we get
$$\align &\sum_{k=1}^n\f{p_k}k-283452.35
\\<&n(\log n+\log\log n-2)-m(\log m+\log\log m-2)
\\&+\log n+\log\log n-1-\f n{\log n}-\f n{\log^2n}+\f{2m}{\log m}
\\&+(\log\log n-1.8)\l(\f n{\log n}+\f{4n}{\log^2 n}+\f{\sqrt n}{\log^2164}-2961.17\r).
\endalign$$
By (2.2),
$$p_n-n>n\l(\log n+\log\log n-2+\f{\log\log n-2.25}{\log n}\r).$$
Hence
$$\align&\sum_{k=1}^n\f{p_k}k-(p_n-n)
\\<&283452.35-m(\log m+\log\log m-2)
\\&+\log n+\log\log n-1-\f n{\log n}-\f n{\log^2n}+\f{2m}{\log m}
\\&+(2.25-1.8)\f n{\log n}+(\log\log n-1.8)\l(\f{4n}{\log^2 n}+\f{\sqrt n}{\log^2164}-2961.17\r)
\\&=-\f {0.55n}{\log n}-\f {8.2n}{\log^2n}+\f{4n\log\log n}{\log^2n}+\f{\sqrt n(\log\log n-1.8)}{\log^2164}
\\&+\log n-2960.17\log\log n-m(\log m+\log\log m-2)+\f{2m}{\log m}
\\&+283451.35+1.8\times2961.17
\\<&-\f {0.55n}{\log n}-\f {8.2n}{\log^2n}+\f{4n\log\log n}{\log^2n}+\f{\sqrt n(\log\log n-1.8)}{\log^2164}
\\&+\log n-2960.17\log\log n+8992.62.
\endalign$$
Since $n\gs 49583$,
$$\f {0.55n}{\log n}+\f {8.2n}{\log^2n}-\f{4n\log\log n}{\log^2n}-\f{\sqrt n(\log\log n-1.8)}{\log^2164}
-\log n+2960.17\log\log n>8992.62$$
and hence $p_n>n+\sum_{k=1}^np_k/k$ as desired.

 (1.2) and (1.3) can be proved in a similar way.  We are done. \qed

\enddocument